\documentclass[conference]{IEEEtran}

\usepackage{gensymb}
\usepackage{comment}
\usepackage{color}
\usepackage{cite}
\usepackage{amsmath, amstext, amssymb}
\usepackage{array}
\usepackage{stfloats}
\usepackage{url}
\usepackage{graphicx}

\usepackage{tikz}
\usetikzlibrary{%
  arrows,
  positioning
}

\IEEEoverridecommandlockouts

\usepackage{multicol,lipsum}

\newcommand{\B}{\mathcal{B}}
\newcommand{\E}{\mathcal{E}}
\newcommand{\G}{\mathcal{G}}
\newcommand{\V}{\mathcal{V}}

\begin{document}

\title{Fast and Robust Determination of Power System Emergency Control Actions}

\author{Sidhant Misra \hspace{0.2in}   Line Roald  \hspace{0.2in}  Marc Vuffray  \hspace{0.2in} Michael Chertkov

\thanks{\hspace*{-15pt} Los Alamos National Laboratory, Center for Nonlinear Systems and T Division, Los~Alamos, NM 87544, \texttt{\{sidhant,roald,vuffray,chertkov\}@lanl.gov}}}

\maketitle

\begin{abstract}
This paper outlines an optimization framework for choosing fast and reliable control actions in a transmission grid emergency situation. We consider contractual load shedding and generation re-dispatch as exemplary emergency actions. To achieve computational efficiency and scalability, this novel formulation of the robust corrective action optimization is stated in terms of voltages and currents, as opposed to standard formulation in terms of power flows and voltages. The current-voltage formulation is natural for expressing voltage and thermal transmission line constraints, and its major advantage is in the linearity of the power flow equations. The complexity of the current-voltage formulation, which is mainly related to the transformation of voltages and currents into active and reactive power injections, can be largely avoided by stating the cost function directly in terms of currents and using linearized or robust estimates for load and generation constraints. The paper considers five different optimization problem formulations, including the full current-voltage formulation, as well as two convex and two linear approximation to the problem. In a case study on an illustrative case study for the IEEE RTS96 system, we demonstrate pros and cons of the different formulations based on simulations for both normal and contingency operation.
\end{abstract}

\begin{IEEEkeywords}
Emergency control, current-voltage formulation
\end{IEEEkeywords}

\IEEEpeerreviewmaketitle

\section{Introduction}
In many parts of the world, increasing loads and fluctuations from renewable energy, public resistance to transmission expansion, and higher focus on economic efficiency is forcing system operators to operate their systems closer to the limits. As a result, many system operators enforce a less strict version of the N-1 criterion, where intermediate post-contingency overloads are tolerated for a short period of time until corrective actions can help the operator to bring the system back to a normal operating state \cite{panciatici2014}. 
While this type of corrective N-1 security allows for less costly system operation \cite{chatzivasileiadis2011, roald2016corrective}, it also increases the possible risk of an emergency situation. The planned corrective controls might fail to react properly \cite{karangelos2013}, or might be insufficient if the system is deviating from the planned operating point (due to, e.g., variability in the renewable generation). Furthermore, the system may experience emergency situations that are more severe than the N-1 contingencies considered in the operational planning, leading to dangerous system operation. In such a situation, where one or more system components are overloaded, the system operator is forced to take quick emergency action. The main goal of those actions is to restore normal operation, remove the immediate danger of component damage and prevent cascading events from further deteriorating the system state. The control actions must be determined in a fast and reliable way, while ensuring that the modification to the current operating point is as small as possible to avoid unnecessarily high cost and long restoration times. 

We are looking for an optimization framework for emergency control which is  both (a) sufficiently {\bf accurate}, in terms of representing reality, 
and (c) sufficiently {\bf efficient}, in terms of computational scalability over realistic (thousands node large) systems, to produce reliable results fast in an emergency situation. In this first formulation, we consider changes in the generation dispatch and contractual load shedding as the possible emergency actions.

The scalability, robustness and efficiency of our newly proposed solution is to be achieved through a number of technical tricks. 
Most importantly, we will utilize the current-voltage, thus called IV, representation of the Power Flow (PF) equations (see \cite{12COC} and references therein). 
This representation is advantageous in that it allows us to easily state both voltage constraints and thermal transmission limits, which are typically given as limits on the current flow. Perhaps more importantly, the IV formulation provides us with a linear representation of the network wide PF constraints, allowing us to remove the impediment of the nonlinearity of the PF equations. 
This advantage of the linear IV approach is important for scalability when the computations are performed over large grids/systems where many degrees of freedom (thousands or even tens of thousands) need to be accounted for. These considerations of the computational scalability has lead researchers recently,  notably in \cite{12OCC, castillo2016}, to restate the basic PF formulation in the IV terms. 

When applied to a typical optimal power flow problem as in \cite{12OCC}, one significant complication of the IV formulation is the (non-convex) transformation of the bus currents into active and reactive power injections. This transformation is necessary to accurately determine the cost of generation, which is minimized as part of the objective function, and to accurately represent generation constraints and the load consumption, which are typically given in active and reactive power. 
In our formulation of the emergency control problem, we will use two tricks to circumvent these issues. 
First, since economic considerations are less important in the emergency situation, we state the objective function directly in terms of currents and voltages and choose penalty functions that, e.g., approximate the cost of generation or limits the number of involved generators. 
Second, generation constraints and loads are represented using convex approximations in the current-voltage domain. In particular, we consider two convex approximations, (i) a conservative, robust inner approximation of the IV feasible domain based on the pre-control system state and (ii) an approximation based on a linear Taylor expansion around the pre-control operating point. We further show how these convex representations can be efficiently linearized to enable computational scalability.


The rest of the paper is organized as follows. In Section~\ref{sec:IV_full} we introduce the emergency control problem and provide a current-voltage formulation of the same. In Section~\ref{sec:convex_approx} and Section~\ref{sec:linear_approx} we develop convex and linear formulations of the emergency control problem in Section~\ref{sec:IV_full}. Section~\ref{sec:simulations} provides a comparison of the methods introduce in this paper via numerical computations. We conclude the paper in Section~\ref{sec:conclusion} and discuss directions for future research.


\section{Current-Voltage Formulation of the Optimal Power Flow Problem} \label{sec:IV_full}
We first introduce some notation that we use throughout the paper. We denote by $\B$ the set of buses, $\E$ the set of transmission lines of a power transmission network. Let $\G$ and $\mathcal{L}$ denote the set of generators and loads. We assume that there is exactly one generator and one load per bus. Presence of multiple generators/loads can be handled by aggregation. 
We denote by $v = (v^i)_{i \in \B}$ as the vector of complex voltage phasors, and $i_B = (i_B^j)_{j \in \B}$ as the vector of bus current injections, where by convention an injection of current is represented by a positive value. The current flow on the branches are given by $(i_F^l)_{l \in \mathcal{L}}$. The generator and load currents are denoted by $(i_G^j)_{j \in \G}$ and $(i_L^j)_{j \in \mathcal{L}}$ respectively. By convention, the generator current $i_G$ enters the bus from the generator, and the load current $i_L$ leaves the bus towards the load. 

The problem we would like to solve, stated in terms of voltages and currents (based on rectangular coordinates $v_{re},~v_{im},~i_{B,re},~i_{B,im},~i_{F,re},~i_{F,im}$), active and reactive power is given by:

\begin{subequations} \label{eq:non_convex}
\begin{eqnarray}
&\min_{i_B,v, i_F, p, q}& \sum_{i\in\mathcal{G,L}}c_{p,i}p_i + c_{q,i}q_i \label{cost_of_correction}\\
&\mbox{s.t.}&
i_B = i_G - i_L, \label{bus_aggregation} \\
&&i_{B} = Y_{B}v, \label{buscurrent}\\
&&i_{F} = Y_{F}v, \label{powerflow}\\
&& \sqrt{i_{F,re}^2+i_{F,im}^2}\leq |i_F|^{max}, \label{currentmax}\\
&& |v|^{min}\leq\sqrt{v_{re}^2+v_{im}^2}\leq |v|^{max}, \label{voltagemag}\\
&& p + j q = v i_B^*, \label{PQtoIV} \\ 
&& p^{min} \leq p \leq p^{max}, \label{genP} \\
&& q^{min} \leq q \leq q^{max}, \label{genQ} \\
&& \theta = 0 \label{slack}
\end{eqnarray}
\label{detproblem}
\end{subequations}
The linear cost function in the objective \eqref{cost_of_correction} represents the cost of emergency control, where the considered control actions are currently limited to generation redispatch, demand response and forced load shedding. The parameters $c_{p},~c_q$ represents the cost coefficients of the respective control actions, i.e., the cost of adjusting generation or demand as well as the cost of load shedding. Note that we expect much higher cost for forced load shedding than for adjustments of the generation or demand, and that the cost of active power adjustments are higher than the cost of reactive power adjustments $c_p>>c_q$. 
Eq.~\eqref{bus_aggregation} denotes the aggregation of generator and load currents at a bus. The nodal power balance is given by \eqref{buscurrent}, where $Y_{B}$ represents the complex bus admittance matrix. Eq. \eqref{powerflow} represents the branch currents $i_F$ from the voltages, with $Y_{F}$ representing the branch admittance matrix. 
The current and voltage magnitudes are constrained by \eqref{currentmax} and \eqref{voltagemag}, respectively. The conversion from nodal active and reactive power consumption to nodal voltage and currents are given by \eqref{PQtoIV}. 
Constraints on active and reactive power generation and consumption is given by \eqref{genP}, \eqref{genQ}. Note that both loads and generation is assumed to be controllable, since the emergency situation mandates the use of load shedding.
Eq. \eqref{slack} sets the reference angle at the slack bus to zero. 

The IV formulation has significant advantages because of the  linearity of  power flow equations \eqref{buscurrent}, \eqref{powerflow}. Instead, the non-convexity of the AC power flow is concentrated in the conversion from $(i,v)$ to $(p,q)$ in \eqref{PQtoIV}, and in the voltage magnitude constraints \eqref{voltagemag}.
In this paper, we argue that in the special circumstances related to emergency control, finding fast control actions to restore system stability is a higher priority than a exact representation of system cost. This enables us to avoid the explicit, non-convex transformation from $(i,v)$ to $(p,q)$, and devise liner or convex programs  that can very efficiently solve the emergency control problem in real time. While the solution might not be optimal in the sense that it might require more generation redispatch or higher load shed than necessary, it is guarantees system feasibility and has reliable and fast solve times. 

In the following sections, we derive four different approximate formulations. These formulations differ in the way they employ simplifications, and approximations, and can be broadly categorized into \emph{Convex} and \emph{Linear} approximations. 
Later in Section~\ref{sec:simulations} we compare the strengths and weaknesses of each of these formulations by numerical computations on the IEEE RTS96 test case.


\section{Convex Approximations of the IV Formulation} \label{sec:convex_approx}
As described above, the main non-convexity of the problem is related to the conversion from the $(i,v)$ to the $(p,q)$ domain at the nodes. 
In this section, we describe two different approximations of this conversion that result in convex programs: (i) an approximation based on Taylor expansion and (ii) a robust inner approximation. 

Since both formulations share the approximation of the cost function \eqref{cost_of_correction} and voltage constraints \eqref{voltagemag}, we first described these constraints. We then describe the handling of the conversion from $(p,q)$ to $(i,v)$ and of the constraints on active and reactive power \eqref{PQtoIV} - \eqref{genQ}. 


\subsubsection{Cost Function}
\label{sec:cost_function}
The cost function represents the cost (or impact) of the emergency control actions, such as the amount of load shed or the cost of generation redispatch. However, while it is important to express a preference to reduce the load shedding to a minimum and limiting the generation redispatch, our main goal is the fast determination of the least impact emergency action rather than an accurate representation of the actual cost. 
Therefore, we approximate both the load shed and the generation redispatch using a first order Taylor expansion around a reference point $(v^0,~i^0)$ of the active and reactive power, 
\begin{align}
    p = &~v_{re}i_{G,re} + v_{im}i_{G,im} \nonumber \\
    \approx &~v^0_{re}i_{G,re} + v_{re}i_{G,re}^0 - v^0_{re}i^0_{G,re} \nonumber  \\
    &~+ v^0_{im}i_{G,im} + v_{im}i_{G,im}^0 - v^0_{im}i^0_{G,im}, \\
    q = &~v_{im}i_{G,re} - v_{re}i_{G,im} \nonumber \\
    \approx &~v_{im}^0 i_{G,re} + v_{im}i_{G,re}^0 - v_{im}^0i_{G,re}^0 \nonumber \\
    &~- v^0_{re}i_{G,im} - v_{re}i_{G,im}^0 + v_{re}^0i_{G,im}^0. \label{eq:taylor}
\end{align}
The reference point can be chosen arbitrarily, but can for example be chosen similar to the current, post-contingency emergency situation or the last recorded system state before the contingency took place. The details of choosing the initial reference point is described in Section~\ref{sec:reference_solution}. 
Given these expressions, we approximate the cost of emergency control \eqref{cost_of_correction} as
\begin{align}
    c_{lin} = &\sum_{i \in \mathcal{G,L}}c_p^i p^i + c_q^i q^i  \nonumber \\
     \approx& c_p^i(v^0_{re}i_{G,re} + v_{re}i_{G,re}^0 - v^0_{re}i^0_{G,re} \nonumber \\  
     &+ v^0_{im}i_{G,im} + v_{im}i_{G,im}^0 - v^0_{im}i^0_{G,im} ) \nonumber \\
     +& c_q^i (v_{im}^0 i_{G,re} + v_{im}i_{G,re}^0 - v_{im}^0i_{G,re}^0 \nonumber \\
    &- v^0_{re}i_{G,im} - v_{re}i_{G_im}^0 + v_{re}^0i_{G,im}^0 ).
    \label{objective_approx}
\end{align}

\subsubsection{Convex Approximation of Voltage Magnitude Constraints}
\label{sec:voltage_convex}
\begin{figure}[ht]
    \centering
    \includegraphics[scale=0.35]{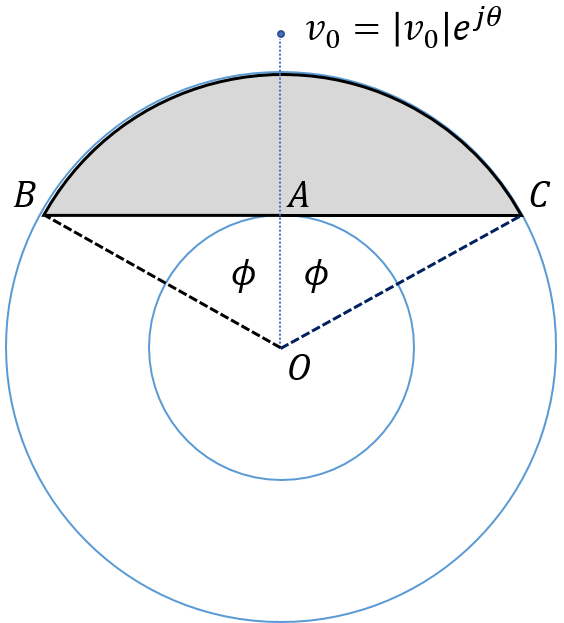}
    \caption{Convex voltage feasibility domain: The voltage domain is obtained by limiting the feasible region to be the largest convex set around the projection of the initial voltage solution $v_0$. The line $BC$ denotes the constraint \eqref{eq:voltage_lower_approx}.}
    \label{fig:voltage_feasibility}
\end{figure}

While the constraint on the current magnitude \eqref{currentmax} is convex (i.e., it forms a circle in the convex plane), 
\begin{equation}
    \sqrt{i_{F,re}^2 + i_{F,im}^2} \leq |i_F|^{max}
\end{equation}
the voltage magnitude constraint \eqref{voltagemag} is non-convex due to the existence of a lower bound $|v|^{min}$,
\begin{equation}
    |v|^{min} \leq \sqrt{v_{re}^2 + v_{im}^2} \leq |v|^{max}. \label{eq:voltage_bounds_nonconvex}
\end{equation}
The non-convex feasible domain for the voltage constraint is illustrated in blue in Fig. \ref{fig:voltage_feasibility}.
We obtain a convex inner approximation of this non-convex domain by constructing a convex feasibility domain (marked in grey in Fig. \ref{fig:voltage_feasibility}, based on the reference voltage solution $v_0$ (as an example, we show the emergency voltage $v^0$ outside of the feasible region).\\ 
The upper bound on voltage magnitude Eq.~\eqref{voltagemag} is a convex quadratic constraint and is kept unaltered. The lower bound is now replaced by a linear constraint as depicted in Figure~\ref{fig:voltage_feasibility}.
To derive the lower bound, observe that the point $A$ in Figure~\ref{fig:voltage_feasibility} is given by $B: \frac{v^0}{|v^0|}|v|^{min}$. The region above the line $BAC$ can then be described using standard tools from coordinate geometry as
\begin{align}
    \left\langle (v_{re}, v_{im}) - \frac{v^0}{|v^0|}|v|^{min}, v^0 \right\rangle \geq 0, \label{eq:voltage_lower_derivation}
\end{align}
where $\langle x,y \rangle$ denotes the inner (dot) product between two vectors $x$ and $y$. Reformulating \eqref{eq:voltage_lower_derivation} and combining with the convex upper bound in \eqref{eq:voltage_bounds_nonconvex}, the inequalities representing the feasible region $S_{\V}$ for the voltage in Figure~\ref{fig:voltage_feasibility} are given by
\begin{subequations} \label{eq:voltage_constraints_approx}
\begin{align}
    \sqrt{v_{re}^2 + v_{im}^2} \leq |v|^{max}, \label{eq:voltage_upper} \\
    v_{re}v_{0,re} + v_{im}v_{0,im} &\geq |v|^{min} |v_0|.\label{eq:voltage_lower_approx}
\end{align}
\end{subequations}
The constraints \eqref{eq:voltage_upper}, \eqref{eq:voltage_lower_approx} are used in all subsequent formulations.

We now describe the two convex formulations of the emergency control computation problem. While both convex formulations use the objective approximation \eqref{objective_approx}, they differ in the representation of the load and generation constraints \eqref{genP}, \eqref{genQ}. Our goal is to represent these constraints on active and reactive power directly in the $(i,v)$ domain to avoid the non-convex conversion to the $(p,q)$ domain. 
Note that in our formulation, both generator and loads are controllable (i.e., we allow for load shedding where the loads can be reduced compared to their initial consumption).

Without loss of generality we present here the approximation to the generator real and reactive power constraints. An identical expression is applicable for the load active and reactive power, so for the following discussions, we omit mentioning the load constraints for brevity. 

\subsection{Convex Formulation using Taylor Expansion for the Power Constraints}
\label{sec:convex_taylor}

In this formulation, we represent the generation and load constraints \eqref{genP}, \eqref{genQ} using first order Taylor expansion \eqref{eq:taylor} around the reference point $(v^0, i_{G}^0)$. With this approximation, the rectangular bounds on active and reactive power generation and consumption in Eqs.~\eqref{genP},\eqref{genQ} become linear:
\begin{align}
    p^{min} &\leq v^0_{re}i_{G,re} + v_{re}i_{G,re}^0 - v^0_{re}i^0_{G,re} \nonumber  \\
    &+ v^0_{im}i_{G,im} + v_{im}i_{G,im}^0 - v^0_{im}i^0_{G,im} \leq p^{max}, \label{eq:real_taylor_approx} \\
    q^{min} &\leq v_{im}^0 i_{G,re} + v_{im}i_{G,re}^0 - v_{im}^0i_{G,re}^0 \nonumber \\
    &- v^0_{re}i_{G,im} - v_{re}i_{G,im}^0 + v_{re}^0i_{G,im}^0 \leq q^{max} \label{eq:reactive_taylor_approx}
\end{align}

The resulting formulation is given below:

\begin{subequations} \label{eq:convex_taylor_formulation}
\begin{align}
    \min \quad &\mbox{cost} \ \eqref{objective_approx} \\
    \mbox{s.t.} \quad &\mbox{network flow equations} \ \eqref{bus_aggregation},\eqref{buscurrent},\eqref{powerflow} \\
    &\mbox{branch current bound} \ \eqref{currentmax} \\
    &\mbox{voltage magnitude bounds} \ \eqref{eq:voltage_constraints_approx} \\
    &\mbox{load/generation power bounds} \ \eqref{eq:real_taylor_approx},\eqref{eq:reactive_taylor_approx}.
\end{align}
\end{subequations}

\subsection{Convex Robust Formulation using Inner Approximation for the Power Constraints}
\label{sec:convex_robust}
In this formulation, we obtain convex inner approximations to \eqref{genP}, \eqref{genQ}. This formulation benefits from feasibility guarantees - a feasible solution to the formulation guarantees a feasible solution to the original non-convex problem in \eqref{eq:non_convex}. However, the robust guarantees requires a certain degree of conservativeness, leading to larger amounts of generation dispatch and possible load shed than is strictly necessary.  


To describe the inner approximation to the active power constraint, we first introduce some notation. Let $\theta \triangleq \arctan\left(\frac{v_{im}^0}{v_{re}^0} \right)$ denote the vector of phase angles of the initial voltage phasor. Let $\phi = \arccos\left( \frac{|v|^{min}}{|v|^{max}} \right)$ denote the collection of angles corresponding to $\angle AOC$ in Figure~\ref{fig:voltage_feasibility}.

\begin{figure}[ht]
    \centering
    \includegraphics[scale=0.3]{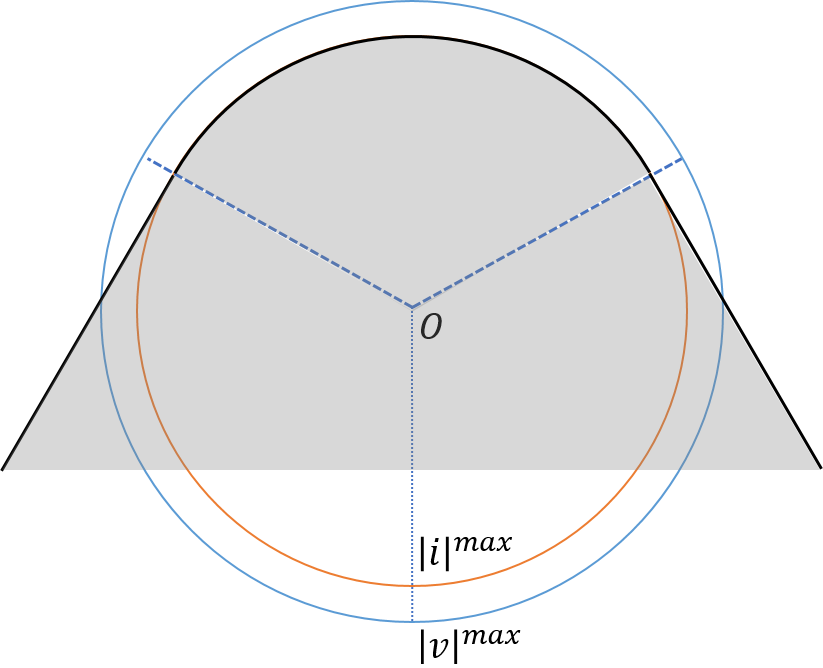}
    \caption{Robust upper current bound: The picture shows the upper robust feasible region given by the constraints \eqref{eq:pmax_robust_hard}. The two straight lines represent the constraints \eqref{eq:pmax_robust_lin1} and \eqref{eq:pmax_robust_lin2}. Although the region defined is convex, it is difficult to directly represent it in the form of intersection of finitely many convex constraints. }
    \label{fig:current_robust_convex_upper}
\end{figure}

Consider the upper bound on the active power given by 
\begin{align}
    v_{re}i_{G,re} + v_{im}i_{G,im} \leq p^{max}. \label{eq:pmax_explicit_nonconvex}
\end{align}
The expression for active power can be thought of as an inner product between the vectors $v_{vec} =  (v_{re}, v_{im})$ and $i_{vec} = (i_{re}, i_{im})$. Using this view, we can re-write the generator active power constraint as 
\begin{align}
    \langle v_{vec}, i_{vec} \rangle \leq p^{max}. \label{eq:inner_product}
\end{align}
Our strategy is to find constraints on the current $i_G$ such that the constraint \eqref{eq:inner_product} is \emph{robustly feasible} for all values of voltage in the voltage feasibility domain $S_{\V}$. From the inner product interpretation, we can identify the worst case $v$ for a given $i_G$. There are three main cases:
\begin{itemize}
    \item $\angle i_G \in [\theta-\phi, \theta+\phi]$: In this case, it is easy to see that the worst-case voltage for constraint \eqref{eq:inner_product} is given by $v = |v|^{max}e^{j \angle i_G}$.
    \item $\angle i_G \in (\theta - \pi, \theta - \phi]$: The worst case voltage is $v = |v|^{max} e^{j(\theta - \phi)}$.
    \item $\angle i_G \in (\theta + \phi, \theta + \pi]$: The worst case voltage is $v = |v|^{max} e^{j(\theta + \phi)}$.
\end{itemize}
The above cases leads to the following convex region for the current $i_G$ described by the following set of constraints, see Figure~\ref{fig:current_robust_convex_upper} for a pictorial representation. 
\begin{subequations} \label{eq:pmax_robust_hard}
\begin{align}
    &i_{G,re}^2 + i_{G,im}^2 \leq \left(\frac{p^{max}}{|v|^{max}}\right)^2, \quad \angle i_G \in [\theta-\phi,\theta+\phi], \label{eq:conditional_current_mag} \\
    &i_{G,re}cos(\theta-\phi) + i_{G,im}sin(\theta-\phi) \leq \frac{p^{max}}{|v|^{max}}, \label{eq:pmax_robust_lin1} \\
    &i_{G,re}cos(\theta+\phi) + i_{G,im}sin(\theta+\phi) \leq \frac{p^{max}}{|v|^{max}}. \label{eq:pmax_robust_lin2}
\end{align}
\end{subequations}
The above feasible region, although simple and convex is difficult to represent directly for passing onto a convex solver primarily due to the conditional constraint \eqref{eq:conditional_current_mag}. Instead it is possible to simply use the unconditional constraint 
\begin{align}
    i_{G,re}^2 + i_{G,im}^2 \leq \left(\frac{p^{max}}{|v|^{max}}\right)^2, \label{eq:unconditional_current_mag}
\end{align}
since it dominates the other two constraints, see Figure~\ref{fig:current_robust_convex_upper}.

\begin{figure}[ht]
    \centering
    \includegraphics[scale=0.3]{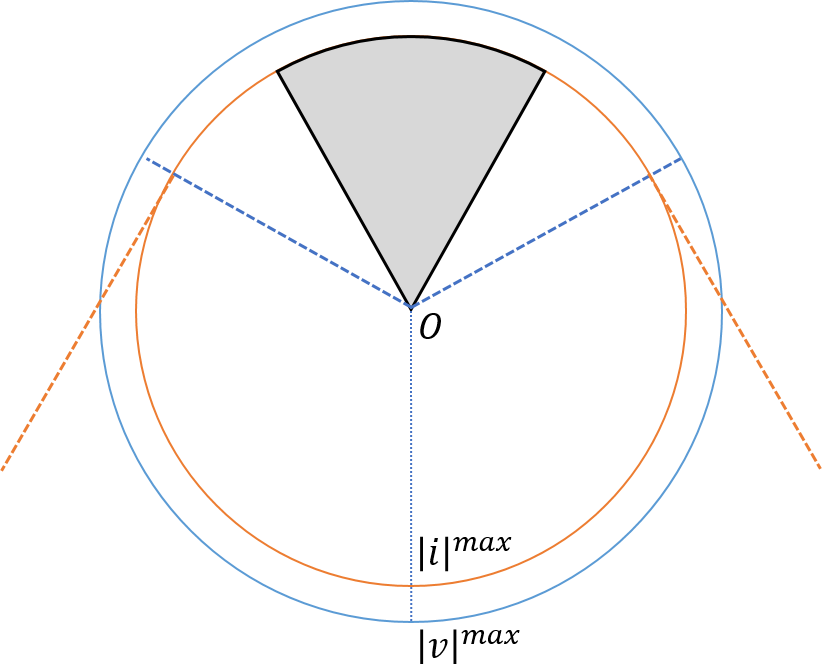}
    \caption{Robust current feasible domain: The figure shows the final current feasible domain obtained by combining the constraints \eqref{eq:pmax_robust_hard} with the constraints \eqref{eq:pmin_robust}.}
    \label{fig:current_robust_convex_combined}
\end{figure}

We proceed similarly to obtain robust current feasibility domain corresponding to the lower bound on the active power. This can be written as 
\begin{align}
        \langle v_{vec}, i_{vec} \rangle \geq p^{min}. \label{eq:inner_product_lower}
\end{align}
To identify the worst case voltage, we again consider the following cases:
\begin{itemize}
    \item $\angle i_G \in [\theta-\pi,\theta]$: The worst case voltage is $v = |v|^{max} e^{j (\theta+\phi)}$.
    \item $\angle i_G \in [\theta,\theta+\pi]$: The worst case voltage is $v = |v|^{max} e^{j (\theta-\phi)}$.
\end{itemize}
Using the above, we get the following robust feasibly constraints for $i_G$:
\begin{subequations} \label{eq:pmin_robust}
\begin{align}
    i_{G,re} cos(\theta-\phi) + i_{G,im}sin(\theta-\phi) &\geq \frac{p^{min}}{|v|^{max}} \\
    i_{G,re} cos(\theta+\phi) + i_{G,im}sin(\theta+\phi) &\geq \frac{p^{min}}{|v|^{max}}.
\end{align}
\end{subequations}
The resulting robust feasible region for the current is depicted in Figure~\ref{fig:current_robust_convex_combined}

We can follow an identical procedure to obtain the robust current feasible domains for the upper and lower bound on the reactive power. Instead of having to repeat the above calculations for the reactive power, we observe that the corresponding constraints can be written in an identical inner product form
\begin{align}
    q^{min} \leq \langle v_{vec}, i_{vec} \rangle \leq q^{max},
\end{align}
where the vector $i_{vec}$ is now modified to read $i_{vec} = (-i_{G,im}, i_{G,re})$. The robust feasible constraints for the reactive power limits can now be written by analogy given by 
\begin{subequations} \label{eq:robust_reactive}
\begin{align}
    &i_{G,im}^2 + i_{G,re}^2 \leq \left(\frac{q^{max}}{|v|^{max}}\right)^2, \label{eq:current_mag_reactive} \\
    &i_{G,re}sin(\theta-\phi) - i_{G,im}cos(\theta-\phi) &\geq q^{min}/|v|^{max} \label{eq:qmin_robust_1} \\
    &i_{G,re}sin(\theta+\phi) - i_{G,im}cos(\theta+\phi) &\geq q^{min}/|v|^{max} \label{eq:qmin_robust_2}.
\end{align}
\end{subequations}

The resulting formulation is given below:
\begin{subequations} \label{eq:convex_robust_formulation}
\begin{align}
    \min \quad &\mbox{cost} \ \eqref{objective_approx} \\
    \mbox{s.t.} \quad &\mbox{network flow equations} \ \eqref{bus_aggregation},\eqref{buscurrent},\eqref{powerflow} \\
    &\mbox{branch current bound} \ \eqref{currentmax} \\
    &\mbox{voltage magnitude bounds} \ \eqref{eq:voltage_constraints_approx} \\
    &\mbox{load/generation power bounds} \ \eqref{eq:pmax_robust_hard},\eqref{eq:robust_reactive}.
\end{align}
\end{subequations}

\section{Linear Approximations of the IV Formulation} \label{sec:linear_approx}
In this section, we derive two fully linear formulations of the emergency control computation problem. Similar to Section~\ref{sec:convex_approx}, these can be classified into: (i) linear approximation using first order Taylor expansion and (ii) linear robust inner approximation. The goal is to obtain linear programs that trade-off accuracy for fast and reliable solutions using modern linear solvers.

To obtain the linear formulations, we need to find linear approximations to the cost \eqref{cost_of_correction}, and the non-linear constraints \eqref{currentmax}, \eqref{voltagemag}, \eqref{genP} and \eqref{genQ}. Both formulations presented here use the linear approximation to the cost function in \eqref{objective_approx}, as well as the linear approximation to the voltage magnitude constraint \eqref{voltagemag} and the branch current magnitude constraint \eqref{currentmax}, which are described in the following sections. 

\subsubsection{Linear Approximation of Branch Current Magnitude Constraints}
\begin{figure}[ht]
    \centering
    \includegraphics[scale=0.35]{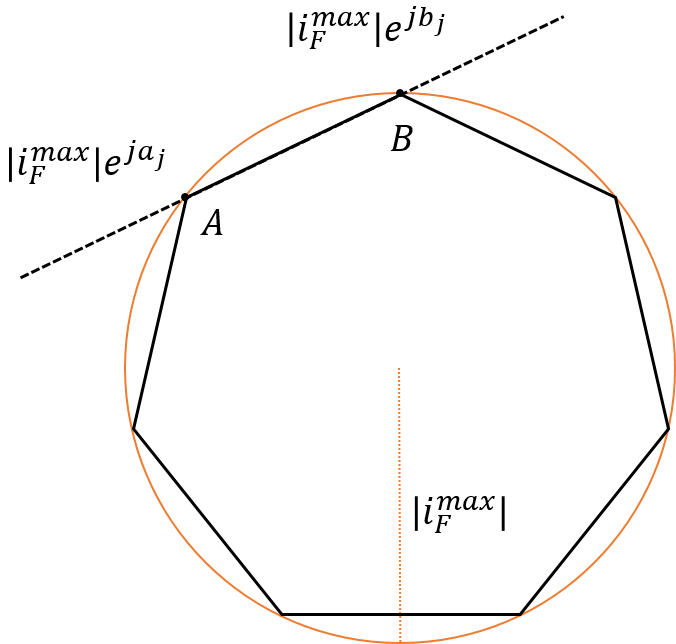}
    \caption{A constraint constituting the linear inner approximation of current magnitude constraint obtained by polyhedral inner approximation to \eqref{currentmax}. The quantities $a_j$ and $b_j$ denote the phase angle of two successive points in the $m_i$ points based uniform discretization of the circle. The inequality in \eqref{eq:current_bound_linear} is denoted by the line $AB$.}
    \label{fig:current_magnitude_inner_linear}
\end{figure}

We derive polyhedral inner approximations to the branch current magnitude constraint \eqref{currentmax} using inscribed regular polygons. The approximation is based on $m_i$ linear inequalities. For $j = 1,2,\ldots, m_i$, let $a_j = (j-1)\pi/m_i$ and $b_j = j\pi/m_i$. Then the $j^{th}$ linear inequality is given by
\begin{align} \label{eq:current_bound_linear}
  i_{F,re} cos(a_j+b_j) + i_{F,im} sin(a_j+b_j) \leq i^{max}cos(a_j-b_j).
\end{align}
The above constraints are shown in Figure~\ref{fig:current_magnitude_inner_linear}. 

\subsubsection{Linear Approximation of Voltage Magnitude Constraints}
\begin{figure}[ht]
    \centering
    \includegraphics[scale=0.35]{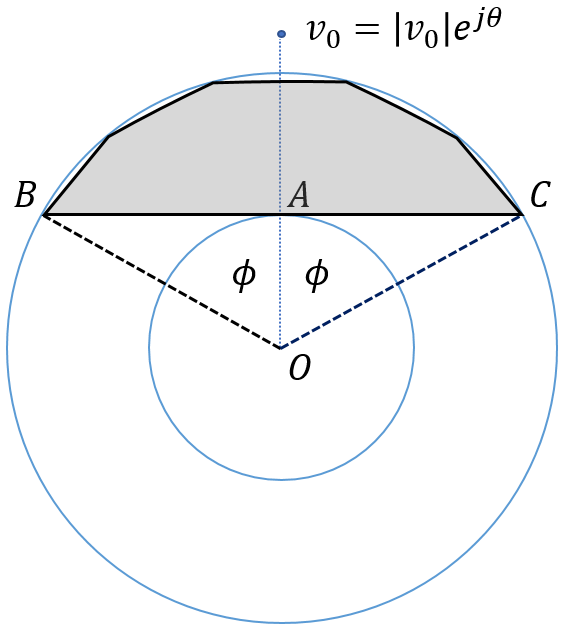}
    \caption{Linear voltage feasibility domain: The linear voltage feasibility domain above is obtained by using a polyhedral approximation to the convex voltage feasibility domain in Figure~\ref{fig:voltage_feasibility} using the constraint in \eqref{eq:voltage_bound_linear}.}
    \label{fig:voltage_feasibility_linear}
\end{figure}

For the voltage magnitude constraints \eqref{voltagemag}, we construct a polyhedral inner approximation to the convex region derived in Section~\ref{sec:voltage_convex}. We observe that we only need to replace the region between $\angle v \in [\theta-\phi,\theta+\phi]$ by a polyhedral inner approximation since the rest of the phase angles are eliminated by the lower bound in \eqref{eq:voltage_lower_approx}. For $j = 1,2,\ldots, m_v$, let 
\begin{align*}
    a_{1,j} &= 0.5(\theta+(j-1)\phi/m_v), \\
    b_{1,j} &= 0.5(\theta+j\phi/m_v), \\
    a_{2,j} &= 0.5(\theta-(j-1)\phi/m_v), \\
    b_{2,j} &= 0.5(\theta-j\phi/m_v).
\end{align*}
Then the $j^{th}$ pair of constraints is given by
\begin{subequations} \label{eq:voltage_bound_linear}
\begin{align}
    v_{re}cos(a_{1,j}+b_{1,j}) &+ v_{im} sin(a_{1,j}+b_{1,j}) \\
    &\leq |v|^{max}cos(a_{1,j}-b_{1,j}),  \\
    v_{re}cos(a_{2,j}+b_{2,j}) &+ v_{im} sin(a_{2,j}+b_{2,j}) \\
    &\leq |v|^{max}cos(a_{2,j}-b_{2,j}).
\end{align}
\end{subequations}
The resuting voltage feasibility domain is given in Figure~\ref{fig:voltage_feasibility_linear}.


\subsection{Linear Formulation using First-Order Taylor Expansion}
The linear formulation in this section is obtained by combining the linear approximations to power using first order Taylor's expansion as described in Section~\ref{sec:convex_taylor} with the voltage and current magnitude approximations described in \eqref{eq:current_bound_linear} and \eqref{eq:voltage_bound_linear}. 
We present the full formulation below:
\begin{subequations} \label{eq:linear_taylor_formulation}
\begin{align}
    \min \quad &\mbox{cost} \ \eqref{objective_approx} \\
    \mbox{s.t.} \quad &\mbox{network flow equations} \ \eqref{bus_aggregation},\eqref{buscurrent},\eqref{powerflow} \\
    &\mbox{branch current bounds} \ \eqref{eq:current_bound_linear} \\
    &\mbox{voltage magnitude bounds} \ \eqref{eq:voltage_bound_linear},\eqref{eq:voltage_lower_approx} \\
    &\mbox{load/generation power bounds} \ \eqref{eq:real_taylor_approx},\eqref{eq:reactive_taylor_approx}.
\end{align}
\end{subequations}

\subsection{Linear Robust Formulation using Inner Approximation}
The linear robust inner approximation uses the same linear inner approximations to the current and voltage magnitude bounds in \eqref{eq:current_bound_linear} and \eqref{eq:voltage_bound_linear} respectively. In addition, we also need to find linear inner approximations to \eqref{eq:conditional_current_mag} or \eqref{eq:unconditional_current_mag}. The linear version no longer suffers from the problem of representing the conditional constraints in \eqref{eq:conditional_current_mag}. As before, using $n_i$ linear pieces, and defining for $j = 1,2,,\ldots,n_i$, and $a_j = (j-1)\pi/n_i$ and $b_j= j\pi/n_i$, we can write the $j^{th}$ linear inequality replacing the conditional constraints in \eqref{eq:conditional_current_mag} as
\begin{align}
      i_{G,re} cos(a_j+b_j) + i_{G,im} sin(a_j+b_j) \leq \frac{p^{max}}{|v|^{max}}cos(a_j-b_j).  \label{eq:robust_current_linear_real}
\end{align}
Similarly for the robust constraint for the reactive power bound in \eqref{eq:current_mag_reactive}, we have the $j^{th}$ inner approximating linear inequality given by 
\begin{align}
    -i_{G,im} cos(a_j+b_j) + i_{G,re} sin(a_j+b_j) \leq \frac{q^{max}}{|v|^{max}}cos(a_j-b_j).   \label{eq:robust_current_linear_reactive}
\end{align}

Another, perhaps more insightful way of obtaining a slightly different form of the robust current constraints \eqref{eq:robust_current_linear_real}, \eqref{eq:robust_current_linear_reactive} is by observing that the original constraint \eqref{eq:pmax_explicit_nonconvex} is linear in $v$ for a given $i_G$. Since in the formulation in this section, the voltage constraints are linear, the values of $v$ that extremize the expression in are the corners of the polytope in Figure~\ref{fig:voltage_feasibility_linear}. With the quantities $a_j, b_j$ described above, the $j^{th}$ corner of the voltage polytope is given by $|v|^{max} e^{j a_j}$. We can now write the $j^{th}$ pair of robust current constraint as 
\begin{align}
    |v|^{max} cos(a_j) i_{G,re} + |v|^{max} sin(a_j) i_{G,im} &\leq p^{max}, \\
    -|v|^{max} cos(a_j) i_{G,im} + |v|^{max} sin(a_j) i_{G,re} &\leq q^{max}.
\end{align}
Unsurprisingly, the expression in the above constraints are very similar to those in \eqref{eq:robust_current_linear_real}, \eqref{eq:robust_current_linear_reactive}.

We present the full formulation below:
\begin{subequations} \label{eq:linear_robust_formulation}
\begin{align}
    \min \quad &\mbox{cost} \ \eqref{objective_approx} \\
    \mbox{s.t.} \quad &\mbox{network flow equations} \ \eqref{bus_aggregation},\eqref{buscurrent},\eqref{powerflow} \\
    &\mbox{branch current bounds} \ \eqref{eq:current_bound_linear} \\
    &\mbox{voltage magnitude bounds} \ \eqref{eq:voltage_bound_linear},\eqref{eq:voltage_lower_approx} \\
    &\mbox{load/gen power upper bounds} \ \eqref{eq:robust_current_linear_real},\eqref{eq:robust_current_linear_reactive} \\
    &\mbox{load/gen power lower bounds} \ \eqref{eq:pmin_robust},\eqref{eq:qmin_robust_1}\eqref{eq:qmin_robust_2}.
\end{align}
\end{subequations}

\section{Choosing the reference point} \label{sec:reference_solution}
The choice of the reference point plays a central role in all four formulations presented in this paper. We propose two different ways of choosing the reference point.

\begin{itemize}
    \item \emph{Pre-contingency system state:} The pre-contingency reference point refers to the state of the system before the contingency that creates an emergency situation has occurred. The
    solution can be easily accessed from the pre-contingency state estimation.
    \item \emph{Post-contingency system state:} The post-contingency reference point refers to the state of the system after the contingency and  prior to deployment of any emergency control actions. The generators and loads that are still online remain at their original value. Since the system is in an emergency state, the post-contingency power flow solution violates one or more line flow and/or bus voltage magnitude constraints. 
\end{itemize}


\section{Case study} \label{sec:simulations}
\begin{figure}[ht]
    \centering
    \includegraphics[scale=0.3]{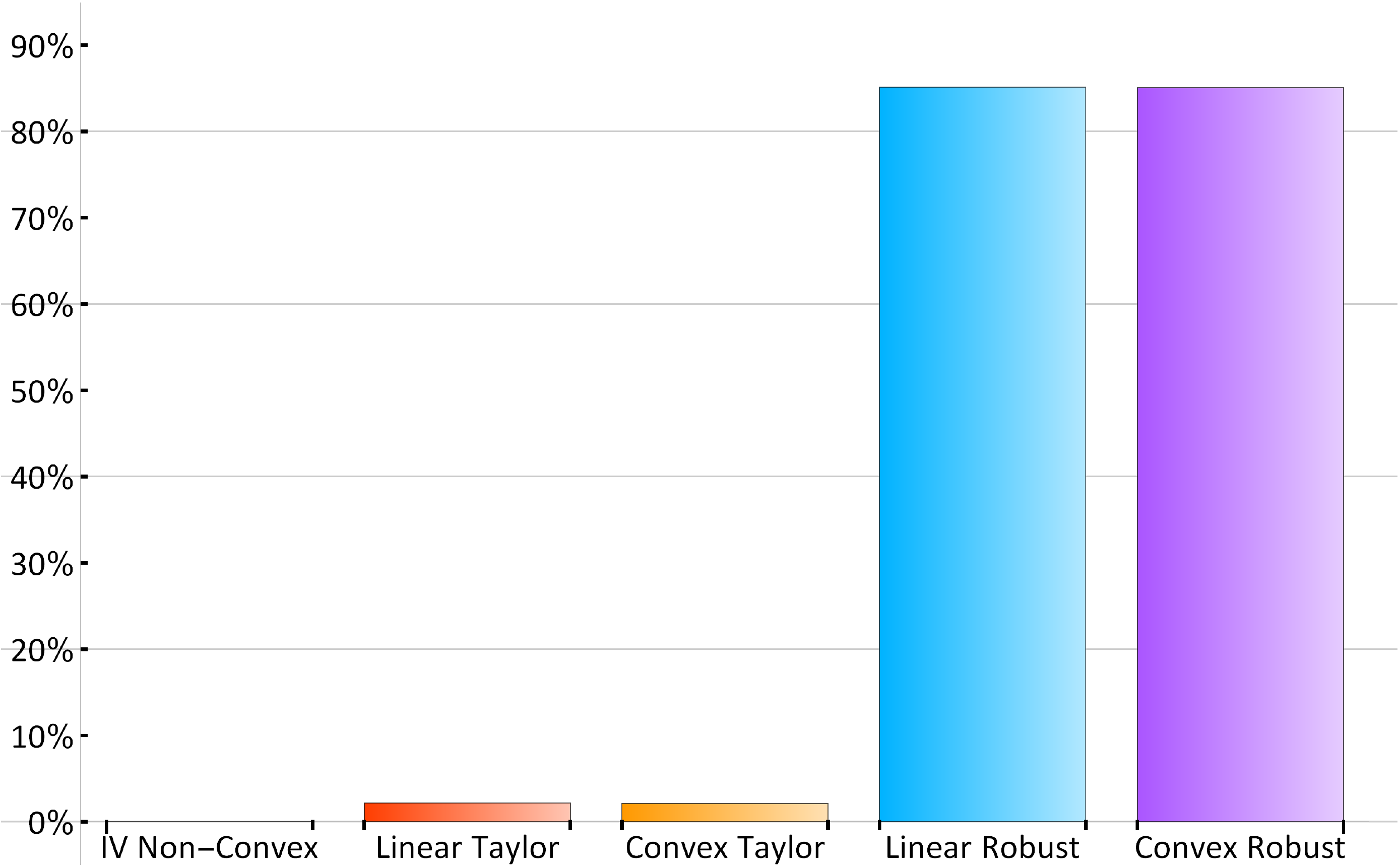}
    \caption{Total load real power shed for each of the four formulations described in this paper as a percentage of the total real power demand. }
    \label{fig:total_load_shed_post}
\end{figure}
\begin{figure}[ht]
    \centering
    \includegraphics[scale=0.3]{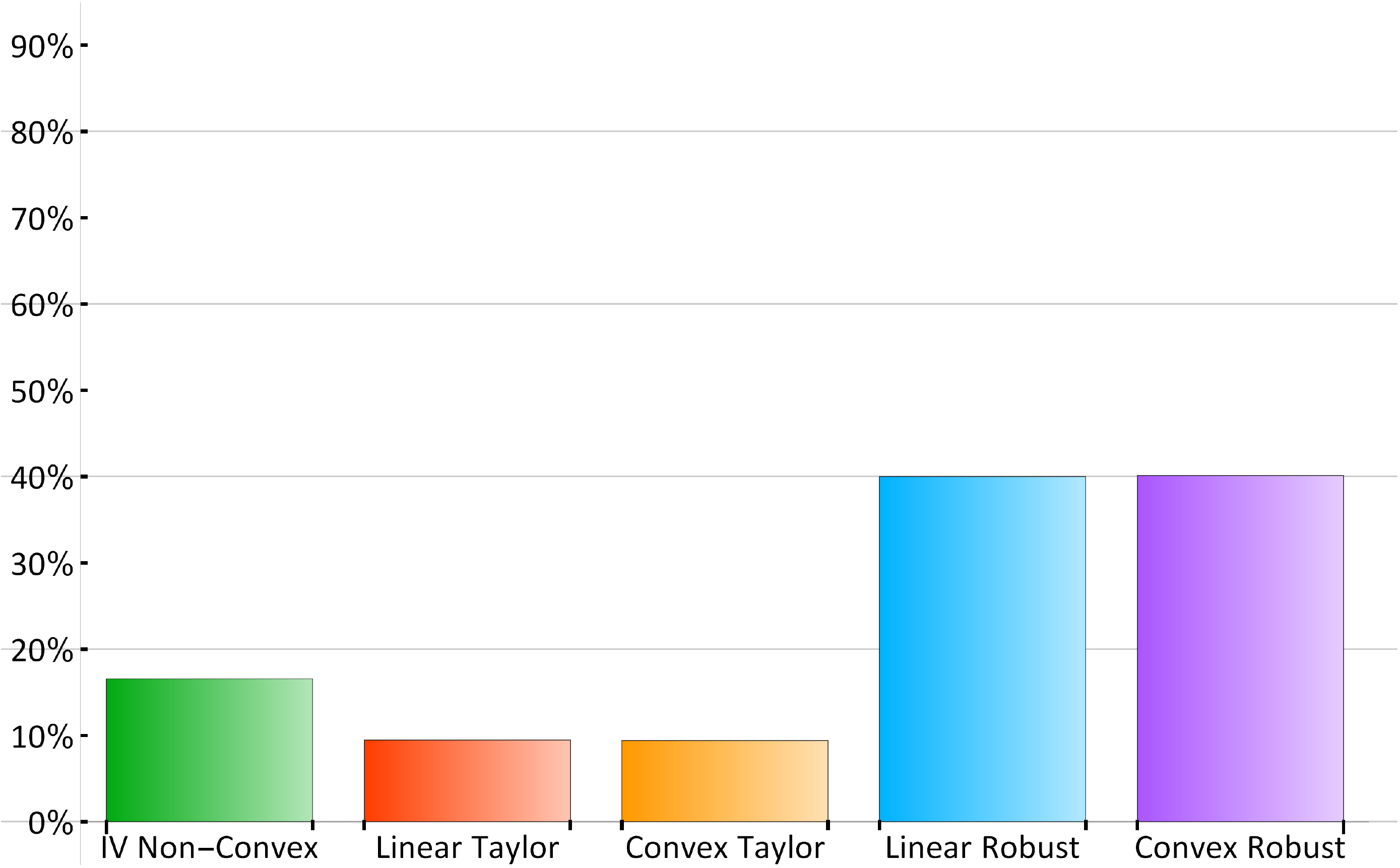}  
    \caption{Total load reactive power shed for each of the four formulations described in this paper as a percentage of the total reactive power demand. }
    \label{fig:total_reactive_load_adjustment}
\end{figure}
\begin{figure}[ht]
    \centering
    \includegraphics[scale=0.3]{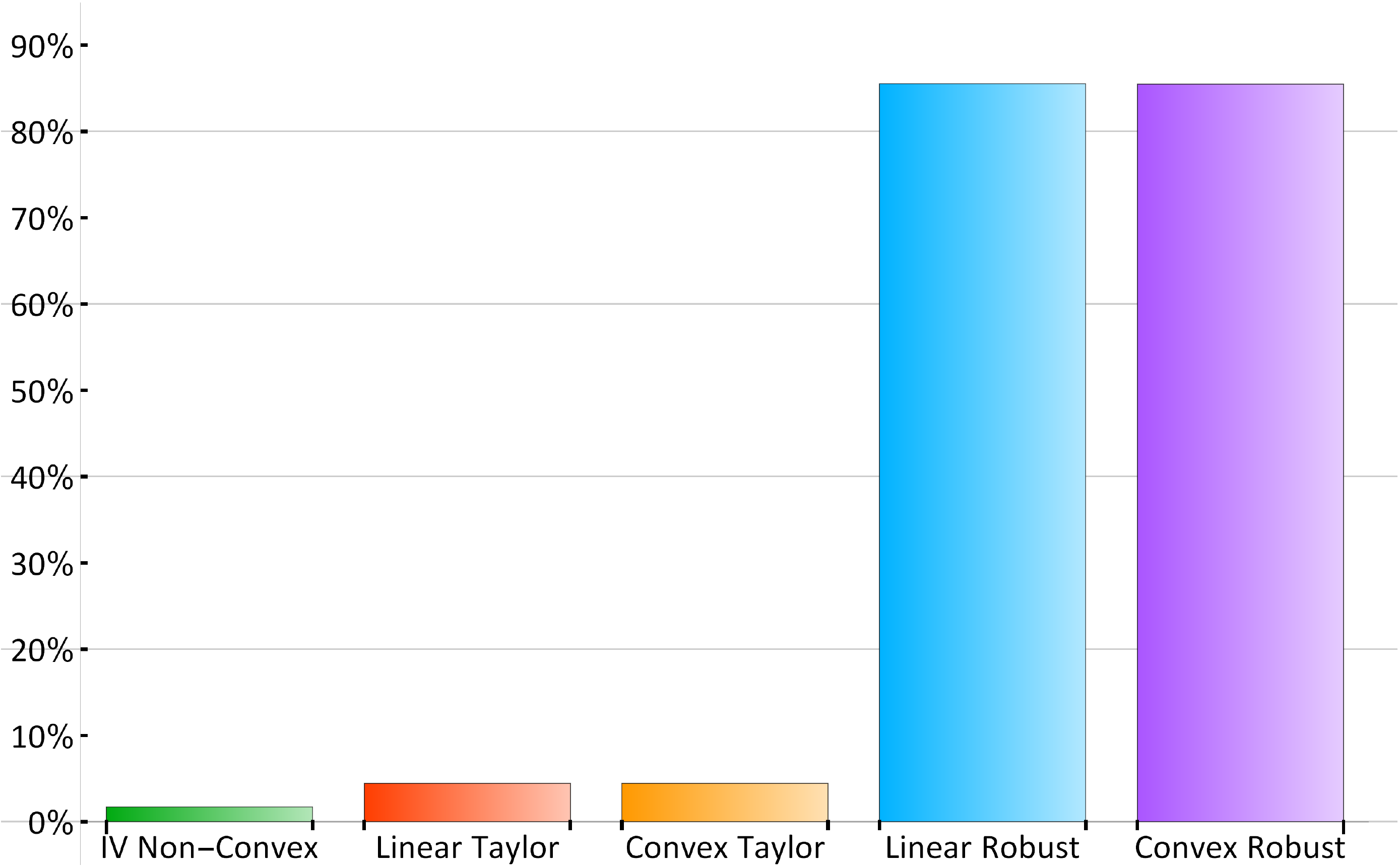}
    \caption{Total active power generation re-dispatch as a percentage of the total active power generation capacity in the system.}
    \label{fig:total_genp_redispatch_post}
\end{figure}
\begin{figure}[ht]
    \centering
    \includegraphics[scale=0.3]{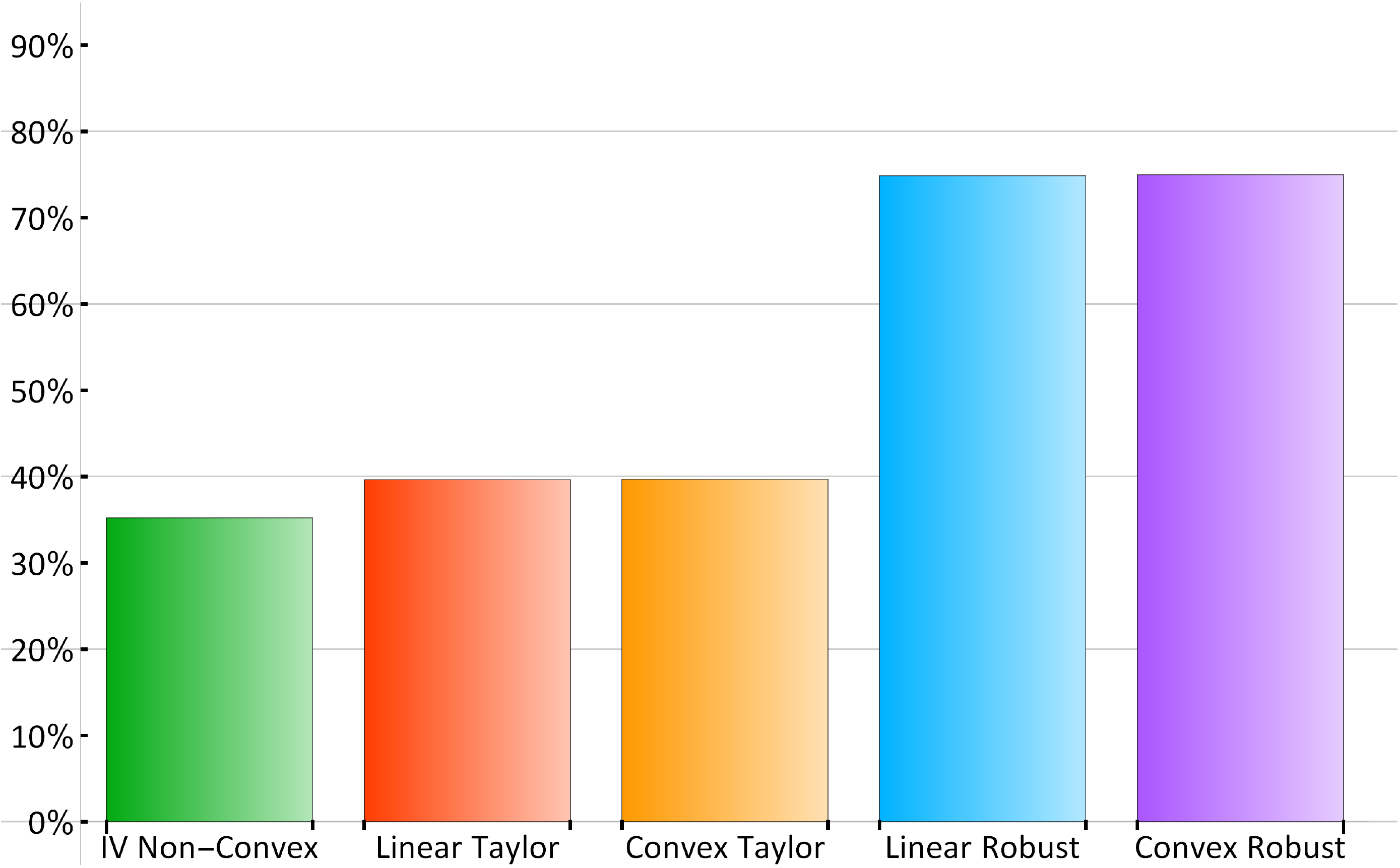}
    \caption{Total reactive power generation re-dispatch as a percentage of the total reactive power generation capacity in the system.}
    \label{fig:total_genq_redispatch_post}
\end{figure}
\begin{figure}[ht]
    \centering
    \includegraphics[scale=0.3]{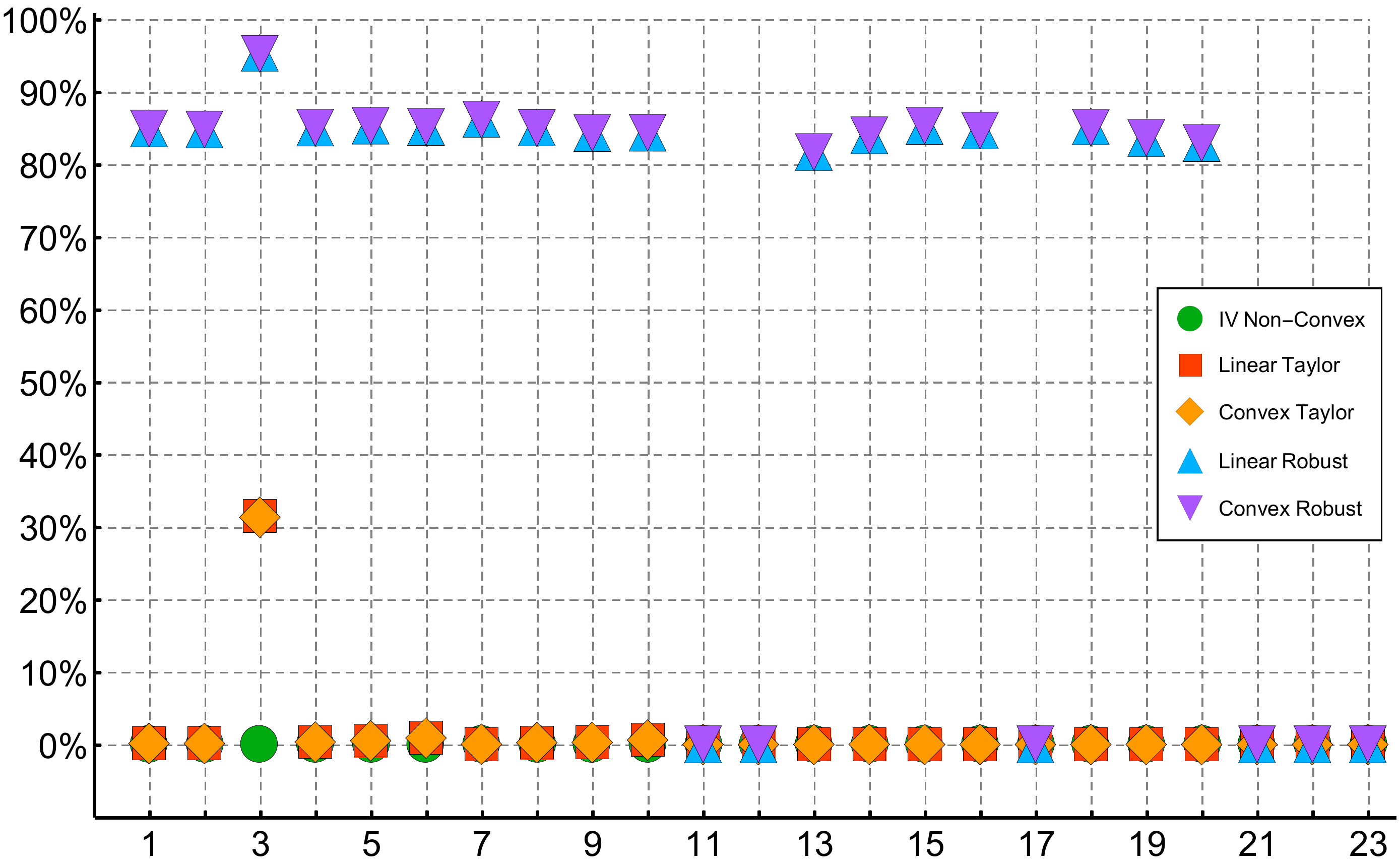}
    \caption{Load real power shed for each bus in the system as a percentage of the real power demand at that bus.}
    \label{fig:load_real_per_bus}
\end{figure}
\begin{figure}[ht]
    \centering
    \includegraphics[scale=0.3]{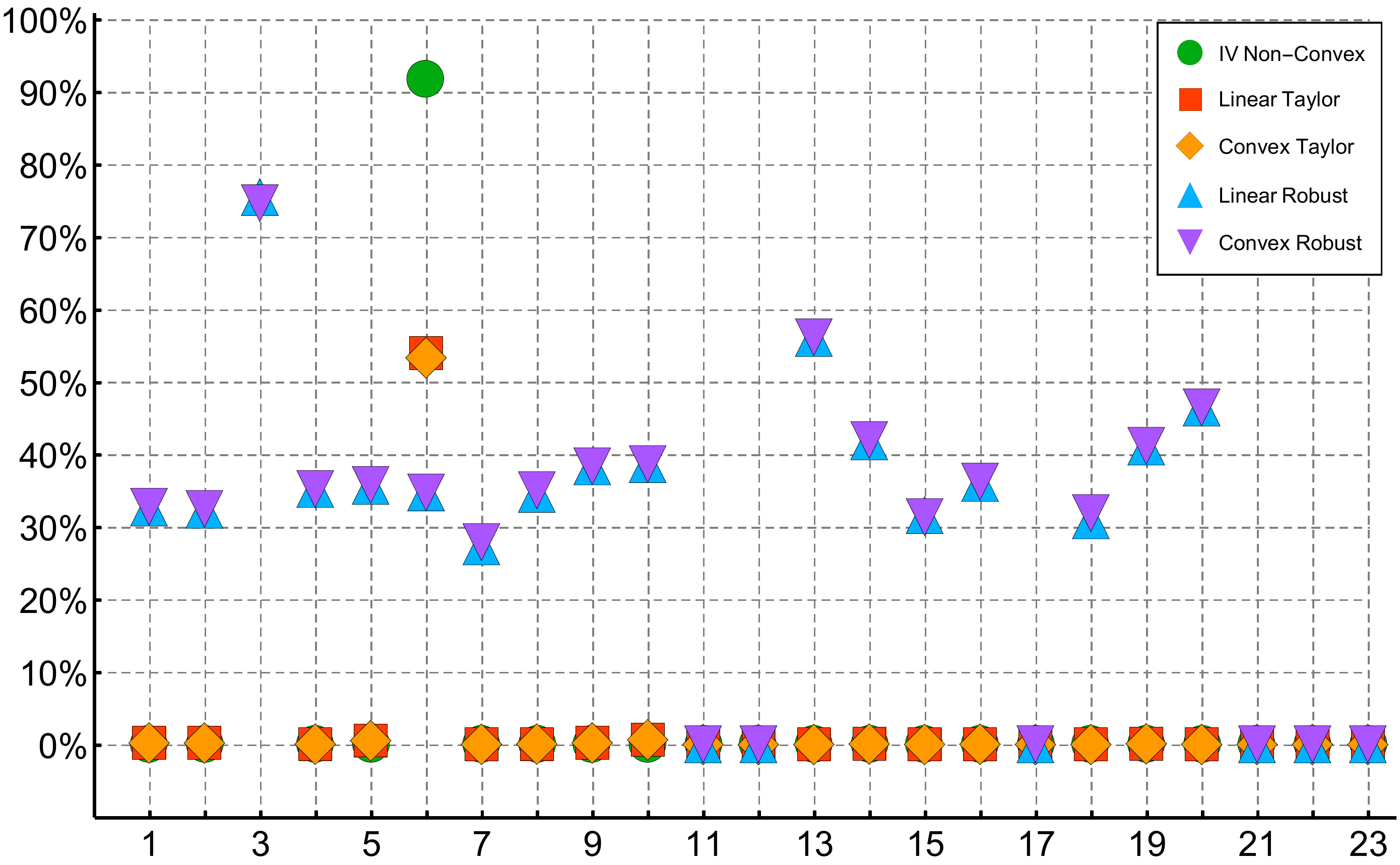}
    \caption{Load reactive power shed for each bus in the system as a percentage of the reactive power demand at that bus.}
    \label{fig:load_reactive_per_bus}
\end{figure}
\begin{figure}[ht]
    \centering
    \includegraphics[scale=0.3]{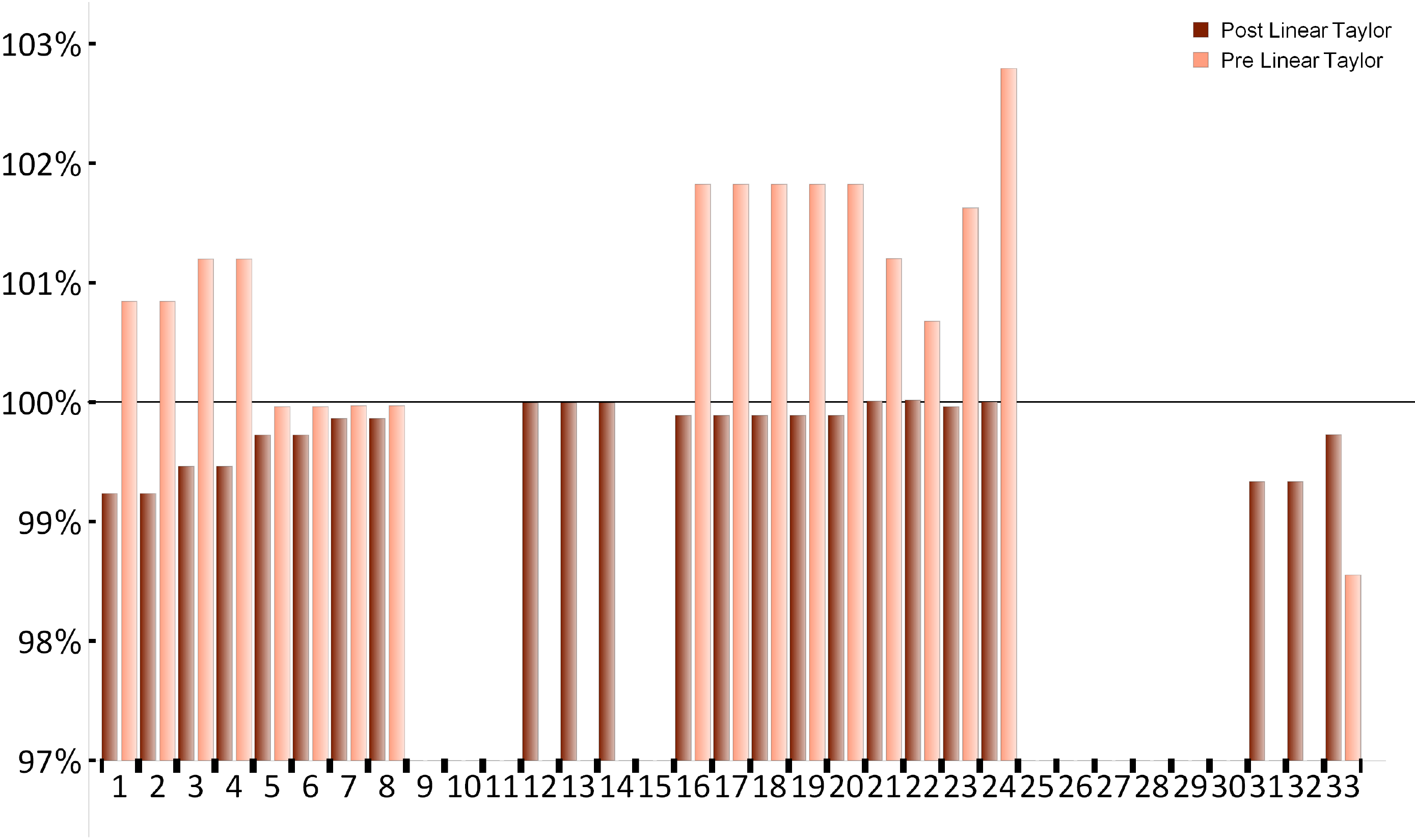}
    \caption{Real power generation for each generator in the system obtained from the linear formulation using Taylor expansion as a percentage of the real power generation limit of that generator. The two cases correspond to reference point being chosen as the pre-contingency and post-contingency solution respectively. The plot has been deliberately magnified to emphasize the generator active power active upper bound violations.}
    \label{fig:pre_post_comparison}
\end{figure}
 In this section, we perform numerical experiments to assess the strengths and weaknesses of each of the formulations described in Section~\ref{sec:convex_approx} and Section~\ref{sec:linear_approx}. The computations are performed on the IEEE RTS96 test system \cite{force1999ieee}. We stress the system by increasing the loads by $15\%$. We then create an emergency state by introducing a bus contingency at bus number $24$, where all the adjacent transmission lines are removed. The data structures from the RTS96 data file are obtained using PowerModels \cite{powermodels} in Julia, and the AC-OPF is solved using the default PowerModels AC-OPF solver in polar coordinates. Note that the standard AC OPF based on polar coordinates does converges to a point of local infeasibility for the post-contingency case (after outage of bus 24), demonstrating the need for more implementations with more reliable convergence. The formulations described in Section~\ref{sec:convex_approx} and Section~\ref{sec:linear_approx} are implemented using JuMP \cite{lubin2015computing}.

We compare the formulations using two metrics: (i) the amount of load shed, and (ii) the amount of generation re-dispatch required to restore feasibility and relieve the emergency state. We also investigate the effect of the choice of the reference point among the options described in Section~\ref{sec:reference_solution}.
Figures~\ref{fig:total_load_shed_post} and \ref{fig:total_reactive_load_adjustment} shows the total load real and reactive power shed in the system for each of the four formulations, when the reference point is chosen to be the post-contingency power flow solution. Figures~\ref{fig:total_genp_redispatch_post} and \ref{fig:total_genq_redispatch_post} shows the same for the total required real and reactive generation power redispatch.

Both robust formulations require significantly higher load shed and generation re-dispatch compared to the formulations using first-order Taylor expansion. However unlike the latter which can suffer from approximation error, the robust formulations are guaranteed to produce feasible solutions. Although the control actions computed from the robust formulations are expensive, the robust formulation can be applied in extreme cases where stabilization and partial rescue of the system is the highest priority. In such scenarios, the robust solution can be employed to quickly restore the system after the partial load shed. Furthermore, the robust formulation can be used as a quick pre-processing step to obtain an initial feasible point for a subsequent optimization step.

To inspect the solutions in more detail, we plot the real and reactive load shed for each bus as a percentage of the original pre-contingency value. These are shown in Figures~\ref{fig:load_real_per_bus} and \ref{fig:load_reactive_per_bus}. We observe that the corresponding convex and linear approximations produce almost identical solutions, implying can be made accurate enough to not introduce substantial error. The stability and computational speed of modern linear programming solvers hence make the linear formulations the appropriate choice for computation of fast emergency control actions. The linear models are also expected to scale extremely well to large realistic transmission system models. In this context, although the IV non-convex formulation outperforms the others in terms of the amount of load shed, the non-convexity of the problem can be expected to lead to poorer scalability as well as lower reliability in providing a solution, since the solvers might either fail to converge or converge to a point of local infeasibility.  Thus, comparing the four convex formulations with respect to the metrics mentioned above, and from a computational stand point, the linear formulation using Taylor expansion described in \eqref{eq:linear_taylor_formulation} outperforms the other formulations.

The reference point plays a crucial role in the derivation of all formulations, and for the best performance, we would like to have a reference point which is close to the final solution. We compare the solutions obtained from the linear formulation using first order Taylor expansion, when the reference point is chosen to be the pre-contingency and the post-contingency system state. To illustrate the difference, we plot the active power generation per bus after the emergency control is computed as shown in Figure~\ref{fig:pre_post_comparison}. 
As can be seen from the figure, several of the normalized active power generation values obtained from choosing pre-contingency solution as the reference are above $100\%$ and thus violate the active power upper limit. On the other hand, no such violation is observed for the generation obtained from choosing post-contingency solution as the reference, making it a better choice.












\section{Conclusion} \label{sec:conclusion}
We considered the problem of computing fast and reliable control actions in power transmission systems in order to relieve the state of emergency caused by contingencies. We propose utilizing the current-voltage formulation where the power flow physics are inherently linear, thus trading off cost for accurate representation of the network flow equations. We formulate four convex/linear approximations to the emergency control computation problem that can be scaled to large systems, and can compute control actions in a fast and reliable manner using modern optimization solvers.  

In future work, we plan to implement and test our formulations on large transmission system models where the number of buses are of the order of a few thousands. Another direction consists in reducing the level of conservativeness of the robust formulations presented in this paper. It might also be useful to explore if algorithms based on hybrid formulations or iterative applications of the formulations in this paper can combine the advantages of the robust (guaranteed feasibility) and Taylor expansion based (lesser load shed and generation re-dispatch) formulations. Furthermore, we would like to investigate constraints that allow for more realistic models of load shed (e.g., methods that preserve the pre-control power factor of the load).

\section*{Acknowledgment}
The work was supported by funding from the U.S. Department of Energy’s Office of Electricity as part of the DOE Grid Modernization Initiative.

\bibliographystyle{IEEEtran}
\bibliography{IREPabstract,NetFlow,StateEst}

\begin{thebibliography}{10}
\providecommand{\url}[1]{#1}
\csname url@samestyle\endcsname
\providecommand{\newblock}{\relax}
\providecommand{\bibinfo}[2]{#2}
\providecommand{\BIBentrySTDinterwordspacing}{\spaceskip=0pt\relax}
\providecommand{\BIBentryALTinterwordstretchfactor}{4}
\providecommand{\BIBentryALTinterwordspacing}{\spaceskip=\fontdimen2\font plus
\BIBentryALTinterwordstretchfactor\fontdimen3\font minus
  \fontdimen4\font\relax}
\providecommand{\BIBforeignlanguage}[2]{{%
\expandafter\ifx\csname l@#1\endcsname\relax
\typeout{** WARNING: IEEEtran.bst: No hyphenation pattern has been}%
\typeout{** loaded for the language `#1'. Using the pattern for}%
\typeout{** the default language instead.}%
\else
\language=\csname l@#1\endcsname
\fi
#2}}
\providecommand{\BIBdecl}{\relax}
\BIBdecl

\bibitem{panciatici2014}
P.~Panciatici, M.~C. Campi, S.~Garatti, S.~H. Low, D.~K. Molzahn, A.~X. Sun,
  and L.~Wehenkel, ``Advanced optimization methods for power systems,'' in
  \emph{Power Systems Computation Conference (PSCC)}, Wroclaw, Poland, Aug.
  2014.

\bibitem{chatzivasileiadis2011}
S.~Chatzivasileiadis, T.~Krause, and G.~Andersson, ``{Flexible AC Transmission
  Systems (FACTS) and Power System Security - A Valuation Framework},'' in
  \emph{{IEEE PES General Meeting}}, Detroit, Michigan, 2011.

\bibitem{roald2016corrective}
L.~Roald, S.~Misra, T.~Krause, and G.~Andersson, ``{Corrective Control to
  Handle Forecast Uncertainty: A Chance Constrained Optimal Power Flow},''
  \emph{IEEE Transactions on Power Systems (in press)}.

\bibitem{karangelos2013}
E.~Karangelos, P.~Panciatici, and L.~Wehenkel, ``Whither probabilistic security
  management for real-time operation of power systems?'' in \emph{2013 IREP
  Symposium on Bulk Power System Dynamics and Control}, Aug 2013, pp. 1--17.

\bibitem{12COC}
\BIBentryALTinterwordspacing
M.~B. Cain, R.~P. O’Neill, and A.~Castillo, ``History of optimal power flow
  and formulations,'' 2012. [Online]. Available:
  \url{https://www.ferc.gov/industries/electric/indus-act/market-planning/opf-papers/acopf-1-history-formulation-testing.pdf}
\BIBentrySTDinterwordspacing

\bibitem{12OCC}
\BIBentryALTinterwordspacing
R.~P. O’Neill, A.~Castillo, and M.~B. Cain, ``The {I}{V} formulation and
  linear approximations of the ac optimal power flow problem,'' 2012. [Online].
  Available:
  \url{https://www.ferc.gov/industries/electric/indus-act/market-planning/opf-papers/acopf-2-iv-linearization.pdf}
\BIBentrySTDinterwordspacing

\bibitem{castillo2016}
A.~Castillo, P.~Lipka, J.~P. Watson, S.~S. Oren, and R.~P. O’Neill, ``A
  successive linear programming approach to solving the iv-acopf,'' \emph{IEEE
  Transactions on Power Systems}, vol.~31, no.~4, pp. 2752--2763, July 2016.

\bibitem{force1999ieee}
R.~T. Force, ``The ieee reliability test system-1996,'' \emph{IEEE Trans. Power
  Syst}, vol.~14, no.~3, pp. 1010--1020, 1999.

\bibitem{powermodels}
\BIBentryALTinterwordspacing
 [Online]. Available: \url{https://github.com/lanl-ansi/PowerModels.jl}
\BIBentrySTDinterwordspacing

\bibitem{lubin2015computing}
M.~Lubin and I.~Dunning, ``Computing in operations research using julia,''
  \emph{INFORMS Journal on Computing}, vol.~27, no.~2, pp. 238--248, 2015.

\end{thebibliography}
\end{document}